\documentclass{article}
\usepackage[pdftex]{hyperref}
\hypersetup{
  colorlinks   = true, 
  urlcolor     = gray, %Color for external hyperlinks
  linkcolor    = red, %Color of internal links
  citecolor   =  blue%Color of citations
}
\usepackage{amsmath,amsfonts,amsthm,amssymb,graphicx,xcolor}
\usepackage[all,2cell,ps]{xy}

\theoremstyle{plain}
\newtheorem{thm}{Theorem}[section]
\newtheorem{lem}[thm]{Lemma}
\newtheorem{prop}[thm]{Proposition}
\newtheorem{cor}[thm]{Corollary}

\theoremstyle{definition}
\newtheorem{defn}[thm]{Definition}

\newtheorem{rem}[thm]{Remark}

\newcommand{\Q}{\mathbb Q}
\newcommand{\C}{\mathbb C}

\DeclareMathOperator{\B}{\mathbf{B}}

\DeclareMathOperator{\supp}{Supp}
\DeclareMathOperator{\mult}{mult}

\newenvironment{pf}{\begin{proof}}{\end{proof}}

\title{On Seshadri constants of varieties with large fundamental group}
\author{\small{Gabriele Di Cerbo}\footnote{Partially supported by the Simons Foundation.} \\ \scriptsize{Columbia University}\\ \footnotesize{\textsf{dicerbo@math.columbia.edu}} \and \small{Luca F. Di Cerbo}\footnote{Partially supported by a grant of the Max Planck Society: ``Complex Hyperbolic Geometry and Toroidal Compactifications'', and by a grant associated to the S. S. Chern position at ICTP.} \\ \scriptsize{International Centre for Theoretical Physics} \\ \footnotesize{\textsf{ldicerbo@ictp.it}}}
\date{}

\begin{document}

\maketitle

%%%%%%%%%%%%%%%%%%%%
\begin{abstract}
Let $X$ be a smooth variety and let $L$ be an ample line bundle on $X$. If $\pi^{alg}_{1}(X)$ is large, we show that the Seshadri constant $\epsilon(p^{*}L)$ can be made arbitrarily large by passing to a finite \'etale cover $p:X'\rightarrow X$. This result answers affirmatively a conjecture of J.-M. Hwang. Moreover, we prove an analogous result when $\pi_{1}(X)$ is large and residually finite. Finally, under the same topological assumptions, we appropriately generalize these results to the case of big and nef line bundles. More precisely, given a big and nef line bundle $L$ on $X$ and a positive number $N>0$, we show that there exists a finite \'etale cover $p: X'\rightarrow X$ such that the Seshadri constant $\epsilon(p^{*}L; x)\geq N$ for any $x\notin p^{-1}\textbf{B}_{+}(L)=\textbf{B}_{+}(p^{*}L)$, where $\textbf{B}_{+}(L)$ is the augmented base locus of $L$.\\ \\ 
\end{abstract}
%%%%%%%%%%%%%%%%%%%%

\tableofcontents\quad\\

%%%%%%%%%%%%%%%%%%%%%%%%%%%%%%%%
\section{Introduction}
%%%%%%%%%%%%%%%%%%%%%%%%%%%%%%%%

Let $X$ be a smooth $n$-dimensional projective variety and let $L$ be an ample line bundle on $X$.
An interesting way of studying the local positivity of the line bundle is by estimating the so-called Seshadri constants of $L$. Recall that, given a nef line bundle $L$ and a point $x\in X$, J.-P. Demailly defines the numbers
\begin{align}\notag
\epsilon(L; x)=\inf_{C\supset x}\frac{L\cdot C}{\mult_{x}(C)},\quad \epsilon(L)=\inf_{x\in X}\epsilon(L; x)
\end{align} 
to be respectively the Seshadri constant of $L$ at $x$ and the global Seshadri constant of $L$, where $C\subset X$ is an irreducible curve and $\mult_{x}(C)$ is the multiplicity of such a curve at $x$. The following proposition explains the connection between separation of jets and Seshadri constants.

\begin{prop}[Demailly, Proposition 6.8. in \cite{Dem90}]\label{Demailly}
Let $X$ be a smooth $n$-dimensional projective variety and let $L$ be a nef line bundle. If $\epsilon(L; x)>n+s$, then the linear system $|K_{X}+L|$ separates $s$-jets at x. Moreover, if $\epsilon(L; x)>2n$ for any $x\in X$, then $K_{X}+L$ is very ample.
\end{prop}

Recall that a line bundle $L$ separates $s$-jets at $x\in X$ if the evaluation map 
\begin{align}\notag
H^{0}(X,L)\rightarrow H^{0}(X,L\otimes \mathcal{O}_{X}/m_{x}^{s+1})
\end{align}
is surjective, where $m_{x}$ is the maximal ideal of the point $x\in X$. In the rest of the paper we will denote by $s(L;x)$ the largest natural number $s$ such that $L$ separates $s$-jets at $x$. One of Demailly's fundamental observations is that $\epsilon(L;x)$ controls the asymptotic of $s(kL;x)$ as a function of $k$, see Proposition \ref{dem2} for a precise statement or the original reference \cite{Dem90} for more details. 

Unfortunately,  the Seshadri constants are in practice very hard to compute or even estimate. The main purpose of this paper is to show that, in the presence of a ``large'' fundamental group, they can be nicely estimated up to a finite cover when $L$ is at least a big and nef line bundle. The following notion for the fundamental group of a variety was introduced by J. Koll\'ar in \cite{KP}.

\begin{defn}\label{Janos}
Let $X$ be a smooth variety. For any positive dimensional irreducible subvariety $Z\subset X$, let us denote by $n_{Z}: \bar{Z}\rightarrow Z$ its normalization. We say that $X$ has large algebraic fundamental group (resp. large fundamental group) if for any such $Z\subset X$ the image $\pi^{alg}_{1}(\bar{Z})\rightarrow\pi^{alg}_{1}(X)$ (resp. $\pi_{1}(\bar{Z})\rightarrow\pi_{1}(X)$ ) is infinite.
\end{defn}

Given Definition \ref{Janos}, we can state our first result.

\begin{thm}\label{second}
Let $X$ be a smooth variety and let $L$ be an ample line bundle on $X$. If $\pi^{alg}_{1}(X)$ is large, given a positive
number $N>0$  there exists a finite \'etale cover $p:X'\rightarrow X$ such that
$\epsilon(p^{*}L)\geq N$.
\end{thm}

Let us observe that Theorem \ref{second} answers affirmatively a conjecture of J.-M. Hwang, see Problem 2.6.2 in \cite{Di Rocco}.\\

In complex differential geometry it is usually more convenient to work with the topological fundamental group, see for example Theorem \ref{Yeung} in Section \ref{preliminaries}. The strategy of the proof of Theorem \ref{second} can be adapted to prove the following.

\begin{thm}\label{first}
Let $X$ be a smooth variety and let $L$ be an ample line bundle on $X$. If $\pi_{1}(X)$ is residually finite and large, given a positive
number $N>0$  there exists a finite \'etale cover $p:X'\rightarrow X$ such that $\epsilon(p^{*}L)\geq N$.
\end{thm}

\begin{rem}
It is possible to show that if $\pi_{1}(X)$ is residually finite and large then $\pi^{alg}_{1}(X)$ has to be large as well. Thus, Theorem \ref{first} is strictly speaking a particular case of Theorem \ref{second}. Nevertheless, we have decided to state Theorem \ref{first} independently of Theorem \ref{second} because in the geometric analysis literature $\pi^{alg}_{1}$ is rarely used. We refer to Theorem \ref{Yeung} below for an example of a classical result in this field.
\end{rem}

Next, Theorems \ref{second} and \ref{first} combined with Proposition \ref{Demailly} have an interesting corollary.

\begin{cor}
Let $X$ be a smooth variety with ample canonical line bundle $K_{X}$. If $\pi_{1}(X)$ is residually finite and large or if $\pi^{alg}_{1}(X)$ is large, there exists a finite \'etale cover $p:X'\rightarrow X$ such that $2K_{X'}$ is very ample.
\end{cor}

Finally, we can state a generalization of Theorems \ref{second} and \ref{first} when $L$ is a big and nef line bundle. 

Recall that given a big line bundle $L$ we denote by $\textbf{B}_{+}(L)$ its augmented base locus, see \cite{Nakamaye} for more details. In Section \ref{big} we prove the following.

\begin{thm}\label{third}
Let $X$ be a smooth variety and let $L$ be a big and nef line bundle on $X$. If $\pi_{1}(X)$ is residually finite and large or if $\pi^{alg}_{1}(X)$ is large, given a positive number $N>0$  there exists a finite \'etale cover $p:X'\rightarrow X$ such that
$\epsilon(p^{*}L; x)\geq N$ for any $x\notin p^{-1}\textbf{B}_{+}(L)=\textbf{B}_{+}(p^{*}L)$.
\end{thm}

Let us observe that when $L$ is ample, i.e. $\textbf{B}_{+}(L)$ is empty, Theorem \ref{third} recovers the statements of Theorems \ref{second} and \ref{first}.\\ \\

\noindent\textbf{Acknowledgements}. The first author would like to express his gratitude to J\'anos Koll\'ar for useful discussions. The second author would like to thank Rita Pardini for generously sharing her knowledge with him and for useful comments and suggestions. He also thanks the Mathematics Department at Notre Dame University and the Max Planck Institute for Mathematics for the nice working environments during the initial stages of this project. Both authors thank the organizers of the conference ``Giornate di Geometria Algebrica e Argomenti Correlati XII'' for getting this project started. \\ \\

%%%%%%%%%%%%%%%%%%%%%%%%%%
\section{Motivations and preliminaries}\label{preliminaries}
%%%%%%%%%%%%%%%%%%%%%%%%%%%%%%%%%

The problem of estimating the Seshadri constants of a ``positive'' line bundle is extensively studied both in the algebraic geometry and geometric analysis literatures.
As shown by Demailly in \cite{Dem90}, these numerical invariants are deeply connected with the theory of singular Hermitian metrics and H\"ormander $L^{2}$-estimates. Therefore they can be studied with techniques coming from geometric analysis. This approach has been successfully explored mainly when the underlying variety  $X$ admits K\"ahler metrics of non-positive sectional curvature. For example, if $(X, \omega_{B})$ is a $n$-dimensional complex hyperbolic manifold equipped with the standard locally symmetric metric Bergman metric whose holomorphic sectional curvature is normalized to be $-1$, J.-M. Hwang and W.-K. To in \cite{Hwang-To} were able to prove that $\epsilon(K_{X}; x)\geq (n+1)\sinh^{2}(i_{x}/2)$ where $i_{x}$ is the injectivity radius of $\omega_{B}$ at $x$.
Similarly, if $(X, \omega)$ is a $n$-dimensional compact  K\"ahler manifold with non-positive sectional curvature, say $-a^{2}\leq R_{\omega}\leq 0$, and $L$ is an ample line bundle on $X$, then S.-K. Yeung in \cite{Yeung1} shows that $\epsilon(L; x)\geq \delta( i_{x}, a, L; n)$ where $\delta$ is an explicitly computable function of the injectivity radius $i_{x}$, the curvature bounds, the curvature of $L$ and the dimension. Moreover, if we fix a normalization for the curvature then $\delta\rightarrow \infty$ as $i_{x}\rightarrow \infty$. This interesting fact then implies the following.

\begin{thm}[Yeung, Theorem 2 in \cite{Yeung1}]\label{Yeung}
Let $(X, \omega)$ be a $n$-dimensional K\"ahler manifold with non-positive sectional curvature and let $L$ be an ample line bundle on $X$. If $\pi_{1}(X)$ is residually finite, given any positive number $N>0$ there exists a finite \'etale cover $p:X^{'}\rightarrow X$ such that $\epsilon(p^{*}L)\geq N$.
\end{thm}

Observe that, by a classical result of H. Wu \cite{Wu}, the universal cover of a compact K\"ahler manifold with non-positive sectional curvature is a Stein manifold. Since Stein manifolds do not contain proper subvarieties, it is easy to show that the varieties considered in Theorem \ref{Yeung} have large fundamental groups, compare with Theorem \ref{first}. This paper started as an attempt of finding a proof of Theorem \ref{Yeung}  which relies on the properties of the fundamental group rather than on curvature assumptions, following a circle of ideas quite common in classical Riemannian geometry.\\

Concerning the organization of the paper, in Section \ref{Preliminaries} we recall the main results in \cite{ELN} and provide the proofs of Theorems \ref{first} and \ref{second}. In Section \ref{big}, we give the details of a proof of Theorem \ref{third}.

\section{Proofs of the main results}\label{Preliminaries}

We start by recalling the connection between Seshadri constants and separation of jets. Let $L$ be a nef line bundle on $X$. For any $k\geq 0$ and given any $x\in X$, let us denote by $s(kL; x)$ the maximal integer such that the linear series $|kL|$ separates $s$-jets at $x$. The following proposition gives a nice lower bound for Seshadri constants in terms of the separation of jets of the linear series $|kL|$. 

\begin{prop}\label{dem2}
Let $X$ be a smooth variety of dimension $n$ and let $L$ be a nef line bundle. For any $k> 0$, we have
\[
\epsilon(L; x)\geq \frac{s(kL; x)}{k}.
\]
\end{prop}

\begin{pf}
Let $C\subset X$ be a reduced irreducible curve containing $x\in X$. By definition for any $k\geq 0$, the linear series $|kL|$ generates $s_{k}$-jets at $x$, where $s_{k}=s(kL; x)$. We can then find a divisor $D_{x}\in |kL|$ which does not contain $C$ such that $\mult_{x}(D_{x})\geq s_{k}$. We then conclude that
\begin{align}\notag
\frac{L\cdot C}{\mult_{x}(C)}\geq \frac{s_{k}}{k}
\end{align}
for any $C\subset X$. Thus, we have
\[
\epsilon(L; x)\geq \frac{s(kL; x)}{k},
\]
and the proof is complete.
\end{pf}

Let us observe that when $L$ is ample then the asymptotic growth of jet separation actually computes $\epsilon(L; x)$ for any given $x\in X$, for more details see Proposition 6.3  in \cite{Dem90} or Chapter V in \cite{Laz1}.  Nevertheless, for our purposes Proposition \ref{dem2} will be sufficient. 

Next, we need to state a result of Ein, Lazarsfeld and Nakamaye proved in \cite{ELN}. This theorem is a variant of the main theorem in the celebrated work of Demailly \cite{Dem93}.

\begin{thm}[Ein-Lazarsfeld-Nakamaye \cite{ELN}]\label{dem1}
Let $X$ be a smooth variety of dimension $n$ and let $L$ be an ample divisor on $X$ satisfying $L^{n}>(n+s)^{n}$. Let $b$ a non-negative number such that $bL-K_{X}$ is nef. Suppose that $m_{0}$ is a positive integer such that $|m_{0}L|$ is base point free. Then for any point $x\in X$ either
\begin{itemize}
\item[-] $|K_{X}+L|$ separates $s$-jets at $x$, or

\item[-] there exists a codimension $c<n$ subvariety $V$ containing $x$ such that
\[
L^{\dim(V)}\cdot V\leq \left( b+m_{0}(n-c)+\frac{n!}{c!}\right)^{c}\left( n+s\right)^{n}.
\]
\end{itemize}
\end{thm}

Now that we recalled these basic results, we can proceed with the proofs of Theorems \ref{second} and \ref{first}. The first step is the construction of a suitably ``large'' \'etale cover of $X$.

\begin{prop}\label{covering lemma}
Let $L$ be an ample line bundle on a smooth variety $X$.  If $\pi_{1}(X)$ is residually finite and large or if $\pi^{alg}_{1}(X)$ is large, given a positive
number $N>0$  there exists a finite \'etale cover $q:X'\rightarrow X$ such that $(q^{*}L)^{\dim(Z)}\cdot Z\geq N$ for any irreducible subvariety $Z\subseteq X'$.
\end{prop}

\begin{pf}
Assume $\pi^{alg}_{1}(X)$ to be large. Let $\hat{\Gamma}$ be the kernel of the map $\pi_{1}(X)\rightarrow \pi^{alg}_{1}(X)$ and let us denote by $\hat{p}: \hat{X}\rightarrow X$ the algebraic universal cover, where $\hat{X}=\tilde{X}/\hat{\Gamma}$  and where $\tilde{X}$ is the topological universal cover.
Let us observe that $\hat{X}$ is non-compact as $\pi^{alg}_{1}(X)$ is infinite.  Next, let $\omega\in [c_{1}(L)]$ be a smooth K\"ahler metric on $X$. Let us consider a  sequence $\{\Gamma_{k}\}$ of nested finite index normal subgroups in $\Gamma=\pi_{1}(X)$  with $\Gamma_{0}=\Gamma$ and such that $\bigcap^{\infty}_{k=0}\Gamma_{k}=\hat{\Gamma}$. Let $\{X_{k}\}$ be the sequence of associated regular \'etale K\"ahler coverings of $X$, where by construction $X_{k}=\hat{X}/\hat{\Gamma}_{k}$ with $\hat{\Gamma}_{k}=\Gamma_{k}/\hat{\Gamma}$.  For any $k$, let us define the following numerical invariant
\begin{align}\notag
\hat{r}_{k}:=\inf\{ d(z, \hat{\gamma}_{k}z)\quad | \quad z\in \hat{X},\quad\hat{\gamma}_{k}\in \hat{\Gamma}_{k},\quad\hat{\gamma}_{k}\neq 1 \},
\end{align}
where the distance is measured with respect to the induced metric on the algebraic universal cover $\hat{X}$. In other words, the distance function $d$ is induced by the K\"ahler metric $\hat{p}^{*}\omega$ on $\hat{X}$.
\begin{lem}\label{Wallach}
We have that $\lim_{k \to \infty}\hat{r}_{k}=\infty$.
\end{lem}
\begin{pf}
If this is not the case, there exist infinite sequences $\{z_{k}\}\in\hat{X}$ and $\hat{\gamma}_{k}\in \hat{\Gamma}_{k}$ such that $d(z_{k}, \hat{\gamma}_{k} z_{k})\leq 2M$ for some positive constant $M$. Let $D$ be a fundamental domain for $X$ in $\hat{X}$, in other words $\hat{p}: D\rightarrow X$ is injective and $\hat{p}: \bar{D}\rightarrow X$ is surjective. Thus for all $k$, there exists $g_{k}\in \hat{\Gamma}_{0}$ such that $g_{k} z_{k}\in\bar{D}$. Let us define $z'_{k}=g_{k} z_{k}$ and $\hat{\gamma}'_{k}=g_{k}\hat{\gamma}_{k}g^{-1}_{k}$, where $\hat{\gamma}'_{k}\in\hat{\Gamma}_{k}$ since $\hat{\Gamma}_{k}$ is by construction a normal subgroup of $\hat{\Gamma}_{0}$. By compactness of $\bar{D}$, there exists a subsequence $\{z'_{k_{j}}\}$ converging to a point $\bar{z}\in\bar{D}$. Since $d(z'_{k}, \hat{\gamma}'_{k}z'_{k})=d(z_{k}, \hat{\gamma}_{k}z_{k})$, we have that 
\begin{align}\notag
d(\bar{z}, \hat{\gamma}'_{k_{j}}\bar{z})\leq d(\bar{z}, z'_{k_{j}})+2M. 
\end{align}
By construction $d(\bar{z}, z'_{k_{j}})\rightarrow 0$, we then conclude that, up to a subsequence, 
$\hat{\gamma}'_{k_{j}}\bar{z}$ converges to a point $q\in\bar{B}(\bar{z}; 2M+\epsilon)$ for some $\epsilon>0$. This implies that 
\begin{align}\notag
\hat{p}(\bar{z})=\hat{p}(\hat{\gamma}'_{k_{j}}\bar{z})\longrightarrow \hat{p}(q)
\end{align}
and then
\begin{align}\notag
(\hat{\gamma}'_{k_{j}}\cdot \hat{\gamma}) q=\hat{\gamma}'_{k_{j}}\bar{z}\longrightarrow q.
\end{align}
for some $\hat{\gamma}\in\Gamma/\hat{\Gamma}$. Now the action of $\hat{\Gamma}_{k}$ on $\hat{X}$ is properly discontinuous, we then conclude that $\hat{\gamma}'_{k_{j}}\cdot \hat{\gamma}=\{1\}$ for all $j$ sufficiently large. Thus, we must have $\hat{\gamma}=\{1\}$ which then implies the contradiction $\hat{\gamma}'_{k_{j}}=\{1\}$. The proof is then complete.
\end{pf} 

By using Lemma \ref{Wallach}, we are now ready to conclude the proof. For any $k$, let us denote by $q_{k}: X_{k}\rightarrow X$ the finite \'etale cover associated to $\Gamma_{k}$. 
Moreover, let us denote by $\hat{p}_{k}: \hat{X}\rightarrow X_{k}$ the algebraic universal covering map. By definition of the numerical invariant $\hat{r}_{k}$, the map $\hat{p}_{k}: B(z; \frac{\hat{r}_{k}}{2})\rightarrow \hat{p}_{k}(B(z; \frac{\hat{r}_{k}}{2}))$ is a biholomorphism for any $z\in\hat{X}$, and then by construction an isometry as well. Thus for any $k$, given a subvariety $Z_{k}\subset X_{k}$, we cannot have $Z_{k}\subset \hat{p}_{k}(B(z; \frac{\hat{r}_{k}}{2}))$ for any ball in $B(z; \frac{\hat{r}_{k}}{2})$ in $\hat{X}$. If otherwise, we can then find a copy of $Z_{k}$ inside the algebraic universal cover $\hat{X}$ which contradicts our topological assumptions on $X$. Recall in fact that $\pi^{alg}_{1}(X)$ is large if and only if $\hat{X}$ does not contain any proper subvariety, see Proposition 2.12 in \cite{KP}. Next, let us observe that all of the metrics $q^{*}_{k}\omega$ have uniformly bounded geometry. In fact, they are pull backs via \'etale maps of a fixed smooth K\"ahler metric on the compact manifold $X$.  In particular, we can find a positive number $r=r(L)>0$ and positive constants $C_{1}$, $C_{2}$
such that 
\begin{align}\label{inequality}
C_{1}\omega_{E}\leq \hat{p}^{*}\omega\leq C_{2}\omega_{E}
\end{align}
on any ball $B(z; r)$ inside $\hat{X}$ and where by $\omega_{E}$ we denote the standard Euclidean K\"ahler metric on $\C^{n}$. Moreover, we can arrange (\ref{inequality}) to hold true for any $\hat{q}^{*}_{k}\omega$ as well as for $\omega$ on balls of the same size. Thus, for any $k$, given an irreducible subvariety $Z_{k}\subset X_{k}$ of pure dimension $m$, there exists a positive constant $K_{m}=K(r, C_{1}, C_{2}; m)$ such that for any point $p\in Z_{k}$ we have that $Vol_{m}(B(p; r)\cap Z_{k})\geq K_{m}$. Here the volume of the $m$-dimensional subvariety $Z_{k}$ is computed by the integral of the $m$-th power of $q^{*}_{k}\omega$ over the smooth part of $Z_{k}$. The volume inequality follows from (\ref{inequality}) and the same statement for the euclidian metric on $\C^{n}$, see Remark \ref{euclide}. Recall now that for any $k$, given a point $p\in Z_{k}$, then the subvariety $Z_{k}$ cannot be entirely contained in the ball $B(p; \frac{\hat{r}_{k}}{2})$. Thus, for any $m$-dimensional irreducible subvariety $Z_{k}\subset X_{k}$ we have that $Vol_{m}(Z_{k})\geq\alpha_{m} \hat{r}_{k}$ for some $\alpha_{m}=\alpha(K_{m})>0$. Next, let us observe that by construction 
\begin{align}\notag
(q_{k}^{*}L)^{m}\cdot Z_{k}=Vol_{m}(Z_{k})\geq \alpha_{m}\hat{r}_{k}
\end{align}
for any $k$. By Lemma \ref{Wallach}, we have that $\hat{r}_{k}\rightarrow\infty$ as $k$ goes to infinity. Thus, given $N>0$ there exists a $k_{0}=k_{0}(N)$ such that for $k\geq k_{0}$ any of the coverings $q_{k}: X_{k}\rightarrow X$ satisfies the requirements of the proposition.

The proof under the assumptions that $\pi_{1}(X)$ is residually finite and large is completely analogous. More precisely, let us consider a sequence $\{\Gamma_{k}\}$ of finite index normal subgroups in $\Gamma=\pi_{1}(X)$ with the property that $\bigcap^{\infty}_{k=0}\Gamma_{k}=\{1\}$. For any $k\geq 0$, we can define the numerical invariant
\begin{align}\notag
r_{k}:=\inf\{ d(z, \gamma_{k}z)\quad | \quad z\in \tilde{X},\quad \gamma_{k}\in \Gamma_{k},\quad \gamma_{k}\neq 1\},
\end{align} 
where $\tilde{X}$ is the topological universal cover, which satisfies $r_{k}\rightarrow\infty$ as $k\rightarrow \infty$. The proof then proceeds exactly as before once we recall that $\pi_{1}(X)$ is large if and only if $\tilde{X}$ does not contain any proper subvariety, see again Proposition 2.12 in \cite{KP}. 

\end{pf}

We are now ready to prove Theorem \ref{second}. The proof of Theorem \ref{first} is identical and we leave its details to the interested reader.

\begin{pf}[Proof of Theorem \ref{second}]
 
Let $n$ be the dimension of the smooth variety $X$. Given the ample line bundle $L$ on $X$, let us fix an integer $a=a(L)$ such that $aL-K_{X}$ is ample and $aL-2K_{X}$ is nef. Moreover, let us define $L'=aL-K_{X}$.
Observe that by Anghern-Siu \cite{Ang} we have that for $m_{0}\geq\binom{n+1}{2}$ then $m_{0}L'$ is base point free.  Next, let us define 
\[
N':=\max\left\{ \left( 1+\binom{n+1}{2}(n-c)+\frac{n!}{c!}\right)^{c}\left( n+s\right)^{n}+1 \: | \: c=0,\dots, n\right\}
\]
where $s>0$ is a fixed parameter to be determined later. By Proposition \ref{covering lemma}, we can find a finite regular \'etale cover $q_{s}:X_{s}\rightarrow X$ such that $(q_{s}^{*}L')^{\dim(V)}\cdot V\geq N'$ for any subvariety $V\subseteq X_{s}$. Moreover, we have that by construction 
\begin{align}\notag
q^{*}_{s}(L'-K_{X})=q^{*}_{s}L'-K_{X_{s}} 
\end{align}
is nef and 
\begin{align}\notag
m_{0}q^{*}_{s}L'=q^{*}_{s}(m_{0}L')
\end{align}
is base point free. By Theorem \ref{dem1}, for any $x\in X_{s}$ the linear system $|K_{X_{s}}+q^{*}_{s}L'|$ separates $s$-jets at $x$.
Since 
\begin{align}\notag
K_{X_{s}}+q^{*}_{s}L'=aq^{*}_{s}L
\end{align}
by Proposition \ref{dem2} we conclude that $\epsilon(q^{*}_{s}L)\geq\frac{s}{a}$. Thus, given any $N>0$ we simply find $s$ such that $\frac{s}{a}\geq N$ and the associated cover $q_{s}: X_{s}\rightarrow X$, which we rename to be $p': X'\rightarrow X$, satisfies the requirements of the theorem. The proof is then complete.

\end{pf}

\begin{rem}\label{euclide}
Let $V$ be a pure $m$-dimensional subvariety passing through the origin $0\in \C^{n}$. Let us consider the Euclidean K\"ahler metric $\omega_{E}$ on $\C^{n}$. We then have the basic inequality
\begin{align}\notag
Vol_{m}(B(0; R)\cap V)\geq mult_{0}(V)\cdot V_{m}(R)
\end{align}
where we denote by $mult_{0}V$ the multiplicity of $V$ at the origin and where $V_{m}(R)$ is the volume of the ball of radius $R$ in $\C^{m}$. For more details see for example page 300 in \cite{Laz1}. Thus, using the uniformly bounded geometries of the covers as in Proposition \ref{covering lemma}, it is possible to directly prove Theorems \ref{first} and \ref{second} without dealing with the separation of jets. We have decided to leave the current proof as it seems of interest to explicitly prove a somewhat stronger result.
\end{rem}

\section{The case of big and nef line bundles}\label{big}

In this section, we appropriately generalize the results of Section \ref{Preliminaries} to the case of big and nef divisors. More precisely, we give the details of the proof of Theorem \ref{third}. Crucial to this generalization is the notion of augmented base locus of a big and nef divisor. This locus is a closed subset of the ambient variety which, roughly speaking, coincides with the set of points where the line bundle ``fails'' to be ample. 
 
Thus, let us start by precisely stating some of the important properties of the augmented base locus of a big line bundle.
To this aim, recall that given a big line bundle $L$ its augmented base locus, denoted by $\B_{+}(L)$, can be
written as
\begin{align}\notag
\textbf{B}_{+}(L)=\bigcap_{L=A+E}\supp(E),
\end{align}
where the intersection is taken over all decompositions $L\sim_{\Q}A+E$,
where $A$ is ample and $E$ effective, see Remark 1.3 in \cite{Nakamaye}.
Note that since $X$ is noetherian and $\B_{+}(L)$ is a Zariski closed subset
of $X$, we can find finitely many decompositions $L\sim_{\Q} A_{i}+E_{i}$ such
that $\B_{+}(L)=\cap_{i=1}^{l}\supp(E_{i})$. 

Next, let us observe that the augmented base locus behaves nicely under finite \'etale covers. More precisely, we can state the following.

\begin{lem}\label{Lemma1}
Let $p: X'\longrightarrow X$ be a regular \'etale cover and let $L$ be a
big and nef line bundle on $X$. We then have
$p^{-1}\textbf{B}_{+}(L)=\textbf{B}_{+}(p^{*}L)$.
\end{lem}
\begin{pf}

By the definition, it is clear that $\textbf{B}_{+}(p^{*}L)\subseteq p^{-1}\textbf{B}_{+}(L)$.
Recall that since $L$ is big and nef by a theorem of Nakamaye, see for example Theorem 10.3.5 in \cite{Laz2}, we know that $\textbf{B}_{+}(L)$ coincides with the union of all positive dimensional subvarieties $V\subseteq X$ such that $L^{\dim(V)}\cdot V=0$. By contradiction, let us now assume the existence of a point $x\in p^{-1}\textbf{B}_{+}(L)$ such that $x\notin \textbf{B}_{+}(p^{*}L)$. Let $y:=p(x)\in X$ which by definition is a point in $\textbf{B}_{+}(L)$. Now, Nakamaye's theorem implies the existence of a positive dimensional subvariety $V$ containing $y$, such that $L^{\dim(V)}\cdot V=0$. Thus, let $W$ be an irreducible component of $p^{-1}(V)$ containing $x$, and let us denote by $d=\textrm{deg}(W\rightarrow V)$ be the degree of the restriction of $p$ to $W$. By the projection formula, we have
\[
((p^{*}L)^{\dim(W)}\cdot W)=d(L^{\dim(V)}\cdot V)=0,
\]
and then, again because of Nakamaye's theorem, it follows that $x\in\textbf{B}_{+}(p^{*}L)$. We have therefore reached a contradiction, and the proof is complete.
\end{pf}

\begin{rem}
We remark that the proof of Lemma \ref{Lemma1} proves a slightly stronger statement. More precisely, it gives that 
\[
p^{-1}\textbf{B}_{+}(L)=\textbf{B}_{+}(p^{*}L),
\]
for $L$ big and nef line bundle on $X$, under the weaker assumption that $p: X'\rightarrow X$ is a finite morphism.
\end{rem}

We can now give the details of the proof of Theorem \ref{third}.

\begin{pf}[Proof of Theorem \ref{third}]
Given $L$, let us fix a decomposition $L\sim_{\Q} A_{i}+E_{i}$ with $A_{i}$ ample and $E_{i}$ effective. Let us choose  $\omega_{i}$ to be a smooth K\"ahler metric in the $\Q$-cohomology class $[A_{i}]$.
By Theorem \ref{first}, given an integer $N>0$, we can construct an infinite sequence of finite regular K\"ahler coverings
\begin{align}\notag
... \rightarrow X_{k+1}\rightarrow X_{k}\rightarrow X_{k-1}\rightarrow ...\rightarrow X,
\end{align}
such that $\epsilon(q^{*}_{k}A_{i})\geq N$ for any $k\geq k_{i}(N, A_{i})$, where $k_{i}$ is a positive integer depending on $N$ and on $A_{i}$. Thus, for any $k\geq k_{i}(N, A_{i})$ and $x\notin q^{*}_{k}(E_{i})$ we have
we have
\begin{align}\notag
\epsilon(q^{*}_{k}L; x)=\epsilon(q^{*}_{k}(A_{i}+E_{i}); x)\geq\epsilon(q^{*}_{k}A_{i}; x)\geq N.
\end{align}
Next, let us find finitely many decompositions $L\sim_{\Q} A_{i}+E_{i}$ such
that $\textbf{B}_{+}(L)=\cap_{i=1}^{l}\supp(E_{i})$. For each $i\in\{1, ..., l\}$ let us chose  $k_{i}$ as above and let us define $K:=Max\{k_{1},..., k_{l}\}$. Next, let us observe that because of Lemma \ref{Lemma1}, given any $x\notin q^{-1}_{k}\textbf{B}_{+}(L)$ we can find at least an index $j\in\{1, ..., l\}$ such that $x\notin q^{-1}_{k}(E_{j})$. Thus, given any $k\geq K$ and for $x\notin q^{-1}_{k}\textbf{B}_{+}(L)$, by choosing an index $j$ as above, we compute
\begin{align}\notag
\epsilon(q^{*}_{k}L; x)=\epsilon(q^{*}_{k}(A_{j}+E_{j}); x)\geq\epsilon(q^{*}_{k}A_{j}; x)\geq N.
\end{align}
The proof is then complete.
\end{pf}

%%%%%%%%%%%%%%%%%%%%

\end{document}